\documentclass{article}
\usepackage{ctex}
\usepackage{amsthm}
\usepackage{amsfonts}
\usepackage{graphicx}
\usepackage{subfigure}
\usepackage{epstopdf}
\usepackage{float}
\usepackage{mathrsfs,makecell}
\usepackage{amstext,graphicx,latexsym,amssymb,amscd,amsfonts,amsmath,slashbox}
\usepackage{float}
\usepackage{color,colortbl,xcolor}
\usepackage{ifthen,calc,lastpage,listings}
\usepackage{amstext, graphicx, latexsym, amssymb, amscd, amsfonts}
\usepackage{booktabs,fancyhdr}
\usepackage{multirow,textcomp,yhmath,exscale,relsize}
\usepackage{ifthen}
\usepackage[left=2.5cm,right=2.5cm,top=2.5cm,bottom=2.5cm]{geometry}
\usepackage[numbers,sort&compress]{natbib}
\usepackage{algorithm}
\usepackage{algorithmic}
\newcolumntype{L}[1]{>{\flushleft}p{#1}}
\newcolumntype{C}[1]{>{\centering}p{#1}}

\AtBeginDvi{\input{zhwinfonts}}
\begin{document}
~\\
\thispagestyle{empty}
\begin{center}
{\huge\textbf{Proximal Quasi-Newton Methods for Multiobjective Optimization Problems$^\ast$}}
\end{center}

\vspace{1mm}
{\large\centerline{Jian-Wen PENG$^\dagger $~~~Jie REN}}
\vspace{2mm}
{\small\centerline{(School of Mathematical Sciences, Chongqing Normal University, Chongqing 401331, China)}}
\vspace{2mm}
{\small\centerline{($^\dagger $E-mail:jwpeng168@hotmail.com)}}

\renewcommand{\thefootnote}{\fnsymbol{footnote}}
\footnotetext[1]{Supported by the National Natural Science Foundation of China (11991024),   the Team Project of Innovation Leading Talent in Chongqing (CQYC20210309536),  the Chongqing University Innovation Research Group Project (CXQT20014) and  the Basic and Advanced Research Project of
Chongqing (CQYC20210309536, cstc2021jcyj-msxmX0300).}
\footnotetext[2]{Corresponding author.}

~\\
{\textbf{ Abstract~~}  \vspace{1mm}

We introduce some new proximal quasi-Newton methods for  unconstrained multiobjective optimization problems (in short, UMOP),  where each objective function   is the sum of a twice continuously differentiable strongly convex function and a proper lower semicontinuous convex but not necessarily differentiable function. We propose proximal BFGS method,  proximal self-scaling BFGS method, and proximal Huang BFGS method for (UMOP) with both line searches and without line searches cases. Under mild assumputions,  we show that each accumulation point of the sequence generated by these   algorithms, if exists, is a Pareto stationary point of the (UMOP).  Moreover, we present their applications in both constrained multiobjective optimization problems and robust multiobjective optimization problems. In particular, for robust multiobjective optimization problems, we show that the subproblems of proximal quasi-Newton algorithms can be regarded as quadratic minimization  problems with quadratic inequality constraints. Numerical experiments are also carried out to verify the effectiveness of the proposed proximal quasi-Newton methods.}

~\\
\textbf{ Key words~~} Multiobjective optimization; Proximal quasi-Newton method; Pareto stationarity; Robust optimization

\vspace{2mm}
\begin{flushleft}
{\zihao{4}\textbf{\Large1\quad Introduction}}
\end{flushleft}

Scalarization approach is one of the most effective methods to solve the multi-objective optimization problem, which transforms the multiobjective optimization problem into a single objective optimization problem (see \cite{G1,GS,Z}). In recent years, the descent method for multiobjective optimization problems has attracted wide attention in the optimization field \cite{FG}.   Fliege  and  Svaiter \cite{FS} proposed  the steepest descent method  for computing a point satisfying the first-order necessary condition of an unconstrained multiobjective optimization problem. Bello Cruz, Lucambio Pérez and Melo \cite{BLM} proposed the projection gradient method for quasiconvex multiobjective optimization problems and showed  that the sequence generated by the algorithm converges to a stationary point. Bonnel, Iusem and Svaiter \cite{BIS} proposed some  proximal point method for vector optimization problems.  Da Cruz Neto, Da Silva, Ferreira, et al. \cite{DDFL}  proposed a subgradient method for quasiconvex multiobjective optimization problems and established the  convergence to Pareto optimal points of the sequences produced by the method. Fliege, Grana Drummond and Svaiter \cite{FGS} introduced the  Newton's method for unconstrained multiobjective optimization problems and showed that the method is   locally superlinear convergent to optimal points. Povalej \cite{P} introduced a quasi-Newton method for unconstrained multiobjective optimization problems and showed that the convergence of this method  is superlinear.

Tanabe, Fukuda  and Yamashita \cite{TFY} introduced some proximal gradient methods  with line searches and without line searches for a unconstrained  multiobjective optimization problem,  where each objective function is the sum of a continuously differentiable convex function and a proper convex  but not necessarily differentiable function. And they also showed that  each accumulation point of these sequence generated by these algorithms, if exists, is Pareto stationary.  They also pointed out in \cite{TFY} that an interesting topic for future research is to propose some new proximal Newton-type algorithms for the  multiobjective optimization problems.

In this paper, we consider the following unconstrained   multiobjective optimization problems (in short, UMOP):
\begin{equation}
\begin{aligned} \label{P}
\quad &\min \quad F(x)  \\
&\begin{array}{r@{\quad}r}
s.t. & x \in {R^n}, \\
\end{array}
\end{aligned}
\end{equation}
where $F:{R^n} \to {(R \cup \{ \infty \} )^m}$ is a vector-valued function with $F: = {({F_1},..,{F_{\rm{m}}})^T}$ and $T$ denotes transpose.  We assume that each ${F_i}:{R^n} \to R \cup \{ \infty \} $ is defined by
\begin{align}
{F_i}(x): = {g_i}(x) + {h_i}(x), ~~~ i = 1,...,m,
\end{align}
where ${g_i}{\rm{:}}{R^n} \to R$ is a twice continuously differentiable strongly convex function, ${h_i}:{R^n} \to R \cup \{ \infty \} $ is proper convex and lower semicontinuous but not necessarily differentiable. It is worthy noting that if $h_i(x) \equiv 0$ for all $x\in {R^n}$ and $i=1, 2, ...,m$, then (1) reduces the multiobjective optimization problems studied in \cite{FG, FS, BLM, BIS, DDFL, FGS, P}.

In order to solve the above (UMOP) (i.e., problem (1) with (2)), we propose some new  proximal quasi-Newton methods with line searches and without line searches based on both the proximal point method and the quasi-Newton method. Here,  the quasi-Newton method is used for the twice continuously differentiable strongly convex function ${g_i}$, and the proximal point method is used for the proper convex but not necessarily differentiable function ${h_i}$.

The main contents of this paper are as follows: In Section 2, we give some notations and some concepts about Pareto optimality and Pareto stationarity. In Section 3, we propose some new proximal quasi-Newton methods, which are called to be   proximal BFGS method,   proximal self-scaling BFGS method, and   proximal Huang BFGS method,  with line searches and without line searches for the (UMOP). In these new proximal quasi-Newton methods, we approximate the Hessian matrices of $g_i$ by using the well known BFGS method, self-scaling BFGS method, and the Huang BFGS method. We prove the global convergence of the proposed algorithms in Section 4. In Section 5, we apply the proposed algorithms to both   constrained multiobjective optimization problems and some robust multiobjective optimization problems. Finally, in Section 6, we verify the effectiveness of the proposed algorithms through numerical experiments to solve robust multiobjective optimization problems.

\vspace{2mm}
\begin{flushleft}
{\zihao{4} \textbf{\Large2\quad Preliminaries}}
\end{flushleft}

For the convenience and brevity of the following discussion, some notations are given in this section, and relevant definitions and lemmas are reviewed.

Let $R$ denote the set of real numbers and $N$ denote the set of positive integers. The Euclidean norm in ${R^n}$ will be denoted by $\left\|  \cdot  \right\|$. We define the relationship $ \le ( < )$ in ${R^m}$ as $u \le v (u < v)$ if and only if ${u_i} \le {v_i} ({u_i} < {v_i})$ for all $i = 1,...,m.$

We call the twice continuously function $g:{R^n} \to R$ is strongly convex if for all $ x,y \in {R^n}$,
\begin{align}
{(\nabla g(x) - \nabla g(y))^T}(x - y) \ge a{\left\| {x - y} \right\|^2}
\end{align}
for some $a > 0$, where $\nabla f(x)$ denote the gradient of $g$ at $x$ (see \cite{NW}).

It's easy to see that (3) is equivalent to

       $~~~~~~~~~~~~~~~~~~~~~~~~~~~~~~~~~~~~~~~  {\nabla ^2}g(x) \ge aI$, for all $x \in {R^n}$, \\
where ${\nabla ^2}g(x)$ denote the Hessian matrix of $g$ at $x$.

It is clear that strong convexity of $g$ implies both the strict convexity and   convexity of $g$. And it is easy to see that  $g$ is strongly convex, its Hessian matrix ${\nabla ^2}{g}(x)$ is positive definite for all $x \in {R^n}$.

Let $f{\rm{:}}{R^n} \to R \cup \{  + \infty \} $, and let $x\in dom(f):=\{x\in R^n: f(x)< +\infty\}$. Then the directional derivative of $f$ at $x$ in the direction $d\in {R^n}$ is defined to be the limit

\[f'(x;d) = \mathop {lim}\limits_{\alpha  \to {0^ + }} \frac{{f(x + \alpha d) - f(x)}}{\alpha },\]
if it exists (see \cite{R}). It's easy to see that $f'(x;d) = \nabla f{(x)^T}d$ when $f$ is differentiable at $x$.

\textbf{Lemma  1} \cite{R}~~Let $f{\rm{:}}{R^n} \to R \cup \{ \infty \} $ be a convex function  and let $x\in dom(f)$. Then for each    $d\in {R^n}$,   the function $\varphi {\rm{:}}(0, + \infty ) \to R$ defined by $\varphi(\alpha )= \frac{{f(x + \alpha d) - f(x)}}{\alpha }$ is non-decreasing. In particular, it follows that
\[f(x + d) - f(x) \ge \frac{{f(x + \alpha d) - f(x)}}{\alpha }~~{\rm for ~all ~} \alpha  \in (0,1).\]

\textbf{Definition  1} \cite{FG, FS}~~Recall that ${x^*} \in {R^n}$ is a Pareto optimal point for (UMOP), if there is no $x \in {R^n}$ such that $F(x) \le F({x^*})$ and $F(x) \ne F({x^*})$. The set of all Pareto optimal values is called Pareto frontier. Likewise, ${x^*} \in {R^n}$ is a weakly Pareto optimal point for (UMOP), if there is no $x \in {R^n}$ such that $F(x) < F({x^*})$.

It's well known that the Pareto optimal point of (UMOP)  is also  a weakly Pareto optimal point of (UMOP), and the converse is not  true in general.

\textbf{Definition 2} \cite{TFY}~ We say that $\bar x \in {R^n}$ is Pareto stationary (or critical) of (UMOP), if and only if
\[\mathop {\max }\limits_{i = 1,...,m} {F'_i}(\bar x;d) \ge 0~ ~{\rm for ~all ~}~d \in {R^n}.\]

It is worthy to noting that Definition 2 generalizes the corresponding ones in \cite{FS} and the following important results hold true.

\textbf{Lemma  2} \cite{TFY}~~(1) If $x \in {R^n}$ is a weakly Pareto optimal point of (UMOP), then $x$ is Pareto stationary.

                          (2) Let every component ${F_i}$ of $F$ be convex. If $x \in {R^n}$ is a Pareto stationary point of (UMOP), then $x$ is also a weakly Pareto optimal point of (UMOP).

                          (3) Let every component ${F_i}$ of $F$ be strictly convex. If $x \in {R^n}$ is a Pareto stationary point of (UMOP), then $x$ is also a Pareto optimal point of (UMOP).

The most popular quasi-Newton's method for nonlinear optimization problems is BFGS method which was introduced by Broyden, Fletcher, Goldfarb and Shanno \cite{B2, F, G2, S}. It is a line search method with a descent direction
\[{d^k} =  - {({B^k})^{ - 1}}\nabla f({x^k}), \]
where $f:R^n \rightarrow R$ is the twice continuously differentiable objective function, ${B^k} \in {R^{n \times n}}$ is an approximation matrix to $ {\nabla}^2 f(x_k)$, which is updated at every iteration as follows:
\begin{align}
{B^{k + 1}} = {B^k} - \frac{{{B^k}{s^k}{{({s^k})}^T}{B^k}}}{{{{({s^k})}^T}{B^k}{s^k}}} + \frac{{{y^k}{{({y^k})}^T}}}{{{{({s^k})}^T}{y^k}}},
\end{align}
where ${s^k} = {x^{k + 1}} - {x^k} $ and ${y^k} = \nabla f({x^{k + 1}}) - \nabla f({x^k})$. As the authors shown in \cite{NW}, ${B^{k + 1}}$ remains positive definite whenever ${B^k}$ is positive definite. The new iterate is
\[{x^{k + 1}} = {x^k} + {\lambda _k}{d^k},\]
where the step length ${\lambda _k} > 0$.

 We  recall the updating formula  of $B^{k}$ related to the self-scaling BFGS (in short, SS-BFGS) method \cite{SY}  as follows:
\begin{align}
{B^{k + 1}} = \frac{{{{({s^k})}^T}{y^k}}}{{{{({s^k})}^T}{B^k}{s^k}}}\left( {{B^k} - \frac{{{B^k}{s^k}{{({s^k})}^T}{B^k}}}{{{{({s^k})}^T}{B^k}{s^k}}}} \right) + \frac{{{y^k}{{({y^k})}^T}}}{{{{({s^k})}^T}{y^k}}}.
\end{align}

It is worthy noting that both the BFGS method and the SS-BFGS method will  be successful if the secant equation ${B^{k + 1}}{s^k} = {y^k}$ and the curvature condition ${({s^k})^T}{y^k} > 0$ are satisfied.

We recall the updating formula  of $B^{k}$ related to  the Huang BFGS (in short, H-BFGS)   as follows:
\begin{align}
{B^{k + 1}} = {B^k} - \frac{{{B^k}{s^k}{{({s^k})}^T}{B^k}}}{{{{({s^k})}^T}{B^k}{s^k}}} + \frac{{{{\hat y}^k}{{({{\hat y}^k})}^T}}}{{{{({s^k})}^T}{{\hat y}^k}}}.
\end{align}

where ${\hat y^k} = {y^k} + \frac{{{\theta ^k}}}{{{{({s^k})}^T}{y^k}}}{y^k}$ and ${\theta ^k}  = 6[f({x^k}) - f({x^{k + 1}})] + 3{[\nabla f({x^k}) + \nabla f({x^{k + 1}})]^T}{s^k}$.

It is also worthy noting that the H-BFGS method  will  be successful if the Huang quasi-Newton equation \cite{SY, ZX} ${B^{k + 1}}{s^k} = {\hat y^k}$    and the curvature condition   ${({s^k})^T}{\hat y^k} > 0$  hold true.

Based on the above statements, if $f$ is a strongly convex function, then the matrix ${B^{k + 1}}$ obtained from each of the mentioned updating formulae for approximating the Hessian matrix always preserves positive definiteness. Thus, in this paper, we assume that all $g_i$ in (UMOP) are strongly convex.

\vspace{2mm}
\begin{flushleft}
{\zihao{4}\textbf{\Large3\quad Proximal quasi-Newton methods}}
\end{flushleft}

Throughout the rest of this paper, we always  assume that  the function $g_i$ is the twice continuously differentiable strongly convex, and the function ${h_i}$ is proper convex and lower semicontinuous but not necessarily differentiable for $i=1, 2,...,m$.

In this section, we propose some new  proximal quasi-Newton methods for (UMOP) with line searches and without line searches.

Now, We define the function ${\theta _x}:{R^n} \to R$ by
\begin{align}
{\theta _x}(d): = \mathop {\max }\limits_{i = 1,...,m} \{ \nabla {g_i}{(x)^T}d + \frac{1}{2}{d^T}{B_i}(x)d + {h_i}(x + d) - {h_i}(x)\} ,
\end{align}
where $\nabla {g_i}(x)$ denotes the gradient of ${g_i}$ at $x$, ${B_i}(x)$ is some approximation of ${\nabla ^2}{g_i}(x)$, $i=1,2,...,m.$ By the convexity of ${g_i}$  and ${h_i}$, we get ${\theta _x}$ is convex and ${\theta _x}(0) = 0.$

We show an important property of ${\theta _x}$ as follows:

\textbf{Lemma 3}~~For all $d \in {R^n},$ the following equality holds:
\[{\theta '_x}(0;d) = \mathop {max}\limits_{i = 1,...,m} {F'_i}(x;d).\]

\begin{proof}[\quad\quad\textbf{Proof}]
Since ${\theta _x}(0) = 0$, By the definition of directional derivative, we get
\[{\theta _x}^\prime (0;d) = \mathop {\lim }\limits_{\alpha  \to {0^ + }} \frac{{{\theta _x}(\alpha d)}}{\alpha }.\]

By the definition of ${\theta _x}$, we have
\begin{align}
\mathop {\lim }\limits_{\alpha  \to {0^ + }} \frac{{{\theta _x}(\alpha d)}}{\alpha } = & {\rm{ }}\mathop {\lim }\limits_{\alpha  \to {0^ + }} \mathop {\max }\limits_{i = 1,...,m} \frac{{\nabla {g_i}{{(x)}^T}(\alpha d) + \frac{1}{2}{{(\alpha d)}^T}{B_i}(x)(\alpha d) + {h_i}(x + \alpha d) - {h_i}(x)}}{\alpha }  \notag \\
= & \mathop {\max }\limits_{i = 1,...,m} \mathop {\lim }\limits_{\alpha  \to {0^ + }} \frac{{\nabla {g_i}{{(x)}^T}(\alpha d) + \frac{1}{2}{{(\alpha d)}^T}{B_i}(x)(\alpha d) + {h_i}(x + \alpha d) - {h_i}(x)}}{\alpha } \notag \\
= & \mathop {\max }\limits_{i = 1,...,m} \{ \nabla {g_i}{(x)^T}d + {h'_i}(x;d)\} \notag \\
= & \mathop {\max }\limits_{i = 1,...,m} {F'_i}(x;d), \notag
\end{align}
where the second equality follows from the continuity of the max function and the third one comes from the definition of directional derivative.
\end{proof}

Let $\omega $ be a positive constant. We define ${\varphi _{\omega ,x}}:{R^n} \to R$ as
\[{\varphi _{\omega ,x}}(d): = {\theta _x}(d) + \frac{\omega }{2}{\left\| d \right\|^2},\]
where the function ${\theta _x}$ is defined in (7). Clearly, ${\varphi _{\omega ,x}}$ is strongly convex and ${\varphi _{\omega ,x}}(0) = 0.$ We would like to define the proximal quasi-Newton direction at an iteration $k$  as ${d^k} = {d_\omega }({x^k})$,   where
\begin{align}
{d_\omega }(x): = \mathop {argmin}\limits_{d \in {R^n}} {\varphi _{\omega ,x}}(d).
\end{align}

\textbf{Remark 1}~~(1) Since ${\varphi _{\omega ,x}}$ is strongly convex, (8) has a unique solution ${d_\omega }(x).$

                   (2) Since ${\varphi _{\omega ,x}}(0) = 0$, we have ${\varphi _{\omega ,x}}({d_\omega }(x)) \le 0.$

Let ${\beta _\omega }(x)$ be the optimal value in (8), i.e.,
\begin{align}
{\beta _\omega }(x): = \mathop {\min }\limits_{d \in {R^n}} {\varphi _{\omega ,x}}(d) = {\varphi _{\omega ,x}}({d_\omega }(x)).
\end{align}

The following lemma characterizes the Pareto stationarity of (UMOP) in terms of ${d_\omega }( \cdot )$ and ${\beta _\omega }( \cdot ).$

\textbf{Lemma 4}~~Let ${d_\omega }(x)$ and ${\beta _\omega }(x)$ be defined in (8) and (9), respectively. Then, the following statements hold.

                (1) If $x$ is a Pareto stationary point of (UMOP), then ${d_\omega }(x) = 0$ and ${\beta _\omega }(x) = 0.$ Conversely, if ${d_\omega }(x) = 0$ and ${\beta _\omega }(x) = 0$, then $x$ is a Pareto stationary point of (UMOP).

                (2) If $x$ is not a Pareto stationary point of (UMOP), then ${d_\omega }(x) \ne 0$ and ${\beta _\omega }(x) < 0$. Conversely, if ${d_\omega }(x) \ne 0$ and ${\beta _\omega }(x) < 0$, then $x$ is not a Pareto stationary point of (UMOP).

                (3) The mappings ${d_\omega }( \cdot )$ and ${\beta _\omega }( \cdot )$ are continuous.

\begin{proof}[\quad\quad\textbf{Proof}]
(1)  Let $x$ be Pareto stationary of (UMOP). Suppose, for the purpose of contradiction, that ${d_\omega }(x) \ne 0$ or ${\beta _\omega }(x) < 0.$ From statements (1) and (2) in Remark 1 it follows that ${d_\omega }(x) \ne 0$ if and only if ${\beta _\omega }(x) < 0$. This means that ${d_\omega }(x) \ne 0$ and ${\beta _\omega }(x) < 0.$

Therefore, we obtain
\begin{align}
{\beta _\omega }(x) = {\theta _x}({d_\omega }(x)) + \frac{\omega }{2}{\left\| {{d_\omega }(x)} \right\|^2} < 0.
\end{align}

It follows from  the convexity of ${\theta _x}$,   ${\theta _x}(0) = 0$ and (10)  that for all $\eta  \in (0,1)$,
\begin{align}
{\theta _x}(\eta {d_\omega }(x)) = & {\theta _x}(\eta {d_\omega }(x) + (1 - \eta ) \cdot 0) \notag \\
\le & \eta {\theta _x}({d_\omega }(x)) + (1 - \eta ){\theta _x}(0)  \notag \\
= & \eta {\theta _x}({d_\omega }(x))  \notag \\
< & - \frac{{\eta \omega }}{2}{\left\| {{d_\omega }(x)} \right\|^2}. \notag
\end{align}

Thus, for all $\eta  \in (0,1)$ we have
\[\frac{{{\theta _x}(\eta {d_\omega }(x))}}{\eta } <  - \frac{\omega }{2}{\left\| {{d_\omega }(x)} \right\|^2}.\]

Since ${d_\omega }(x) \ne 0$, letting $\eta  \to {0^ + }$ we obtain
\[{\theta '_x}(0;{d_\omega }(x)) \le  - \frac{\omega }{2}{\left\| {{d_\omega }(x)} \right\|^2} < 0.\]

It then follows from Lemma 3 that
\[\mathop {\max }\limits_{i = 1,...,m} {F'_i}(x;{d_\omega }(x)) < 0,\]
which contradicts the Pareto stationarity of $x$.

 Let us now prove the converse. Then, suppose that ${d_\omega }(x) = 0$ and ${\beta _\omega }(x) = 0$. From the definition of ${\beta _\omega }(x)$ given in (9), we have
\[{\varphi _{\omega ,x}}(d) = {\theta _x}(d) + \frac{\omega }{2}{\left\| d \right\|^2} \ge {\beta _\omega }(x) = 0~ ~{\rm for ~all ~}~d.\]

Let $\eta  \in (0,1),$ we get
\[\frac{{{\theta _x}(\eta d) + \frac{\omega }{2}{{\left\| {\eta d} \right\|}^2}}}{\eta } \ge 0~ ~{\rm for ~all ~}d.\]

Letting $\eta  \to {0^ + }$ and by Lemma 3, we know that for all $d\in R^n$
\[  {\theta '_x}(0;{d}) =  \mathop {\max }\limits_{i = 1,...,m} {F'_i}(x;d) \ge 0,\]
which implies that $x$ is a   Pareto stationary point of (UMOP).

(2) This statement is equivalent to statement (1).

(3) It is easy to see that the function
\[\mathop {\max }\limits_{i = 1,...,m} \{ \nabla {g_i}{(x)^T}d + \frac{1}{2}{d^T}{B_i}(x)d + {h_i}(x + d) - {h_i}(x)\}  + \frac{\omega }{2}{\left\| d \right\|^2}\]
is continuous with respect to $x$ and $d$. Therefore, the optimal value function ${\beta _\omega }( \cdot )$ is also continuous from \cite[Maximum Theorem]{BP}. Moreover, since the optimal set mapping ${d_\omega }( \cdot )$ is unique, ${d_\omega }( \cdot )$ is continuous from \cite[Corollary 8.1]{H}.
\end{proof}

\begin{flushleft}
{\zihao{-4}\textbf{3.1~Some new proximal quasi-Newton methods for (UMOP) with line searches}}
\end{flushleft}

Now, we present some new proximal quasi-Newton methods for (UMOP) with line searches. To compute the step length ${\lambda _k} > 0$, we use an Armijo rule. Let $\tau   \in (0,1)$ be a prespecified constant. The condition to accept ${\lambda _k}$ is given by
\begin{align}
{F_i}({x^k} + {\lambda _k}{d^k}) \le {F_i}({x^k}) + {\lambda _k}\tau {\theta _{{x^k}}}({d^k}),i = 1,...,m.
\end{align}

We begin with ${\lambda _k} = 1$ and while (11) is not satisfied, we update
\[{\lambda _k}: = \zeta {\lambda _k},\]
where $\zeta  \in (0,1)$. The following lemma demonstrates the finiteness of this procedure.

\textbf{Lemma 5}~~Let ${d^k}$ be defined in (8) with $x = x{}^k$ and $\tau  \in (0,1)$, If ${x^k}$ is not Pareto stationary, then there exists some ${\bar \lambda _k} > 0$ such that
\[{F_i}({x^k} + \lambda {d^k}) \le {F_i}({x^k}) + \lambda \tau {\theta _{{x^k}}}({d^k}),i = 1,...,m\]
for any $\lambda  \in (0,{\bar \lambda _k}]$.

\begin{proof}[\quad\quad\textbf{Proof}]
Let $\lambda  \in (0,1].$ Since ${h_i}$ is convex for all $i = 1,...,m,$ we have
\begin{align}
{h_i}({x^k} + \lambda {d^k}) - {h_i}({x^k}) = & {h_i}((1 - \lambda ){x^k} + \lambda ({x^k} + {d^k})) - {h_i}({x^k}) \notag\\
\le & (1 - \lambda ){h_i}({x^k}) + \lambda {h_i}({x^k} + {d^k}) - {h_i}({x^k}) \notag \\
= & \lambda ({h_i}({x^k} + {d^k}) - {h_i}({x^k})). \notag
\end{align}	
		
Therefore, from the second-order Taylor expansion of ${g_i}$ we obtain
\begin{align}
&{g_i}({x^k} + \lambda {d^k}) + {h_i}({x^k} + \lambda {d^k}) \notag \\
\le & {g_i}({x^k}) + \lambda \nabla {g_i}{({x^k})^T}{d^k} + \frac{1}{2}{(\lambda {d^k})^T}{B_i}({x^k})(\lambda {d^k}) + {h_i}({x^k}) + \lambda ({h_i}({x^k} + {d^k}) - {h_i}({x^k})) + o({\lambda ^2}) \notag \\
= & {g_i}({x^k}) + {h_i}({x^k}) + \lambda [\nabla {g_i}{({x^k})^T}{d^k} + \frac{\lambda }{2}{({d^k})^T}{B_i}({x^k})({d^k}) + {h_i}({x^k} + {d^k}) - {h_i}({x^k})] + o({\lambda ^2}) \notag \\
\le & {g_i}({x^k}) + {h_i}({x^k}) + \lambda [\nabla {g_i}{({x^k})^T}{d^k} + \frac{1}{2}{({d^k})^T}{B_i}({x^k})({d^k}) + {h_i}({x^k} + {d^k}) - {h_i}({x^k})] + o({\lambda ^2}) \notag \\
\le & {g_i}({x^k}) + {h_i}({x^k}) + \lambda {\theta _{{x^k}}}({d^k}) + o({\lambda ^2}) \notag \\
= & {g_i}({x^k}) + {h_i}({x^k}) + \lambda \tau {\theta _{{x^k}}}({d^k}) + \lambda \left[ {(1 - \tau ){\theta _{{x^k}}}({d^k}) + \frac{{o({\lambda ^2})}}{\lambda }} \right], \notag
\end{align}
where ${B_i}({x^k})$ is some approximation of ${\nabla ^2}{g_i}({x^k}),i=1,...,m,$ the second inequality follows from the positive definiteness of ${B_i}({x^k})$ and $\lambda  \in (0,1]$, and the third one comes from the definition of ${\theta _x}$. Since ${x^k}$ is not Pareto stationary, we have ${\theta _{{x^k}}}({d^k}) < 0$ from Lemma 4. It follows from $\tau \in (0,1)$ that there exists some ${\bar \lambda _k} > 0$ such that
\[{g_i}({x^k} + \lambda {d^k}) + {h_i}({x^k} + \lambda {d^k}) \le {g_i}({x^k}) + {h_i}({x^k}) + \lambda \tau {\theta _{{x^k}}}({d^k}),i = 1,...,m,\forall \lambda  \in (0,{\bar \lambda _k}].\]
\end{proof}

To simplify the notation we will use $B_i^k$ to denote ${B_i}({x^k})$  for all $i = 1,...,m$ and $k=0, 1,2,...$.

Based on the previous discussions,  we would like to state our new proximal quasi-Newton methods with line searches for (UMOP) as follows:

\textbf{Algorithm 1}

Step 1 Choose $\omega  > 0$, $\tau  \in (0,1)$, $\zeta  \in (0,1)$, ${x^0} \in {R^n}$, symmetric positive definite matrix $B_i^0 \in {R^{n \times n}},i = 1,...,m$ and set $k: = 0$;

Step 2 Compute ${d^k}$ by solving subproblem (8) with $x = {x^k}$;

Step 3 If ${d^k} = 0$, then stop. Otherwise, proceed to the next step;

Step 4 Compute the step length ${\lambda _k} \in (0,1]$ as the maximum of \[{\Lambda _k}: = \{ {\lambda } = {\zeta ^j}|j \in N,{F_i}({x^k} + {\lambda }{d^k}) \le {F_i}({x^k}) + {\lambda }\tau {\theta _{{x^k}}}({d^k}),i = 1,...,m\} ;\]

Step 5 Set ${x^{k + 1}} = {x^k} + {\lambda _k}{d^k}$, update $\{B_i^{k }\}$  by either one of following three formulae
\begin{align}
B_i^{k + 1} = B_i^k - \frac{{B_i^k{s^k}{{({s^k})}^T}B_i^k}}{{{{({s^k})}^T}B_i^k{s^k}}} + \frac{{y_i^k{{(y_i^k)}^T}}}{{{{({s^k})}^T}y_i^k}},
\end{align}

\begin{align}
B_i^{k + 1} = \frac{{{{({s^k})}^T}y_i^k}}{{{{({s^k})}^T}B_i^k{s^k}}}\left( {B_i^k - \frac{{B_i^k{s^k}{{({s^k})}^T}B_i^k}}{{{{({s^k})}^T}B_i^k{s^k}}}} \right) + \frac{{y_i^k{{(y_i^k)}^T}}}{{{{({s^k})}^T}y_i^k}},
\end{align}

or
\begin{align}
B_i^{k + 1} = B_i^k - \frac{{B_i^k{s^k}{{({s^k})}^T}B_i^k}}{{{{({s^k})}^T}B_i^k{s^k}}} + \frac{{\hat y_i^k{{(\hat y_i^k)}^T}}}{{{{({s^k})}^T}\hat y_i^k}},
\end{align}
where ${s^k} = {x^{k + 1}} - {x^k} = {\lambda _k}{d^k}$, $y_i^k = \nabla {g_i}({x^{k + 1}}) - \nabla {g_i}({x^k})$, $\hat y_i^k = y_i^k + \frac{{\theta _i^k}}{{{{({s^k})}^T}y_i^k}}y_i^k$ and $\theta _i^k: = 6[{g_i}({x^k}) - {g_i}({x^{k + 1}})] + 3{[\nabla {g_i}({x^k}) + \nabla {g_i}({x^{k + 1}})]^T}{s^k}$.
Set $k: = k + 1$, and go to Step 2.

It is worthy noting that, when using (12), (13) and (14) to update $\{B_i^{k }\}$, the above proximal quasi-Newnon methods for (UMOP) will  be  called PQNM(BFGS) method with line searches, PQNM(SS-BFGS)  method with line searches and   PQNM(H-BFGS) method with line searches, respectively.

Observe that from Lemma 4, Algorithm  1 stops at Step 3 with a Pareto stationary point or produces an infinite sequence of nonstationary points $\{ {x^k}\} $. If Step 4 is reached in some iteration $k$, it means that in Step 3, ${d^k} \ne 0$, or equivalently, $ {\beta _\omega }({x^k}) < 0.$  Thus, we have ${\theta _{{x^k}}}({d^k}) < 0.$ It follows from the Armijo condition that objective values sequence $\{{F}({x^k}) \}$ is $R_+^m$-decrease, i.e., \[{F}({x^{k+1}} ) < {F}({x^k})  ~{\rm for ~all ~k}.\]

\begin{flushleft}
{\zihao{-4}\textbf{3.2~Some new proximal quasi-Newton methods for (UMOP) without line searches}}
\end{flushleft}

In this section, we assume that  $\nabla {g_i}$ is Lipschitz continuous with constant $L$ for all $i=1, 2,...,m$. Let the step length ${\lambda _k} \equiv 1$ for $k=0, 1, 2,...$. Now we introduce some new proximal quasi-Newton methods for (UMOP) without line searches.

\textbf{Algorithm 2}

Step 1 Choose $\omega  > L/2$, ${x^0} \in {R^n}$, symmetric positive definite matrix $B_i^0 \in {R^{n \times n}},i = 1,...,m$ and set $k: = 0$;

Step 2 Compute ${d^k}$ by solving subproblem (8) with $x = {x^k}$;

Step 3 If ${d^k} = 0$, then stop. Otherwise, proceed to the next step;

Step 4 Set ${x^{k + 1}} = {x^k} + {d^k}$, update $B_i^{k + 1}$ by either one of the formula (12), (13) or (14). Set $k: = k + 1$, and go to Step 2.

  When using (12), (13) and (14) to update $\{B_i^{k }\}$, we call the proximal quasi-Newton methods for (UMOP) in Algorithm 2 to be
PQNM(BFGS) method without line searches, PQNM(SS-BFGS)  method without line searches and   PQNM(H-BFGS) method without line searches, respectively.
And it is easy to see that the algorithm 2 stops at Step 3 with a Pareto stationary point or generates an infinite sequence of nonstationary points $\{ {x^k}\} $. Moreover, as we can see from the proof of Lemma 9, the objective function values also decrease in each iteration, i.e., \[{F_i}({x^k} + {d^k}) < {F_i}({x^k}),i = 1,...,m.\]

\vspace{2mm}
\begin{flushleft}
{\zihao{4}\textbf{\Large4\quad Convergence analysis}}
\end{flushleft}

In this section, we prove that the sequences generated by both Algorithm 1 and Algorithm 2  converge to Pareto stationary points of (UMOP).

\textbf{Lemma 6} (Three points property) \cite[Lemma 3.2]{CT}~~Let $\sigma :{R^n} \to R \cup \{ \infty \} $ be proper convex and define
\[{x^*} = \mathop {\arg \min }\limits_{x \in {R^n}} \left\{ {\sigma (x) + \frac{1}{2}{{\left\| {x - y} \right\|}^2}} \right\}.\]
Then, for all $ z \in {R^n}$, we have
\[\sigma ({x^*}) - \sigma (z) \le  - \frac{1}{2}{\left\| {z - {x^*}} \right\|^2} - \frac{1}{2}{\left\| {y - {x^*}} \right\|^2} + \frac{1}{2}{\left\| {z - y} \right\|^2}.\]

\textbf{Lemma 7}~~Let $\{ {d^k}\}$ be generated by Algorithms 1 or 2. Then, we have
\[{\theta _{{x^k}}}({d^k}) \le  - \omega {\left\| {{d^k}} \right\|^2}~~{\rm for ~all ~}~k,\]
where ${\theta _x}$ in defined in (7).

\begin{proof}[\quad\quad\textbf{Proof}]
Defining $\sigma(d) : = {\theta _{{x^k}}}({d})/\omega $, we can rewrite (8) with $x = {x^k}$ as
\[{d^k} = \mathop {\arg \min }\limits_{d \in {R^n}} \left\{ {\sigma (d) + \frac{1}{2}{{\left\| {d - 0} \right\|}^2}} \right\}.\]

Thus, substituting ${x^*} = {d^k}$ and $y = z = 0$ into Lemma 6, we get
\[\sigma ({d^k}) - \sigma (0) \le  - {\left\| {{d^k}} \right\|^2}.\]

It follows from  ${\theta _{x{}^k}}(0) = 0$ that
\[{\theta _{x{}^k}}({d^k}) \le  - \omega {\left\| {{d^k}} \right\|^2}~~{\rm for ~all ~}~k.\]
\end{proof}

\begin{flushleft}
{\zihao{-4}\textbf{4.1~Convergence of Algorithm 1}}
\end{flushleft}

\textbf{Lemma 8}~~Let $\{ {d^k}\} $ be generated by Algorithm 1 and suppose that $\{ {F_i}({x^k})\} $ is bounded from below for all $i = 1,...,m$. Then, it follows that
\[\mathop {\lim }\limits_{k \to \infty } {\lambda _k}{\left\| {{d^k}} \right\|^2} = 0.\]

\begin{proof}[\quad\quad\textbf{Proof}] It follows from Lemma 7 and step 4 of Algorithm 1   that
\[{F_i}({x^k} + {\lambda _k}{d^k}) \le {F_i}({x^k}) - {\lambda _k}\tau \omega {\left\| {{d^k}} \right\|^2},i = 1,...,m.\]

Adding up the above inequality from $k = 0$ to $k = \hat k$, where $\hat k$ is a positive integer, we obtain
\begin{align}
~{F_i}({x^{\hat k + 1}}) \le {F_i}({x^0}) - \tau \omega \sum\limits_{k = 0}^{\hat k} {{\lambda _k}{{\left\| {{d^k}} \right\|}^2}} .
\end{align}

Since $\{ {F_i}({x^k})\} $ is bounded from below for all $i = 1,...,m$,  there exists ${\hat F_i} \in R$ such that ${\hat F_i} \le {F_i}({x^k})$ for all $i$ and $k$.

It follows from (15) that
\begin{align}
\sum\limits_{k = 0}^{\hat k} {{\lambda _k}{{\left\| {{d^k}} \right\|}^2}}  \le & \frac{1}{{\tau \omega }}({F_i}({x^0}) - {F_i}({x^{\hat k + 1}}))\notag \\
\le & \frac{1}{{\tau \omega }}({F_i}({x^0}) - {\hat F_i}). \notag
\end{align}

Taking $\hat k \to \infty $, we have $\sum\limits_{k = 0}^\infty  {{\lambda _k}{{\left\| {{d^k}} \right\|}^2}}  < \infty $

and hence $\mathop {\lim }\limits_{k \to \infty } {\lambda _k}{\left\| {{d^k}} \right\|^2} = 0$.
\end{proof}

\textbf{Theorem 1} (i)~Suppose that $\{ {F_i}({x^k})\} $ is bounded from below for all $i = 1,...,m$.~Then every accumulation point of the sequence $\{ {x^k}\} $ generated by Algorithm 1, if it exists, is a Pareto stationary point.

 (ii) Moreover, if the level set of ${F}$ in the sense that $\{ x \in {R^n}\mid F(x) \le F({x^0})\} $ is bounded, then $\{ {x^k}\} $ has accumulation points and they are all Pareto stationary.

\begin{proof}[\quad\quad\textbf{Proof}] We now prove the  first statement.

 Let $\bar x$ be an accumulation point of $\{ {x^k}\} $ and let $\{ {x^{{k_j}}}{\rm{\} }}$ be a subsequence converging to $\bar x$. From statement (3) of Lemma 4, we have ${d^{{k_j}}} = {d_\omega }({x^{{k_j}}}) \to {d_\omega }(\bar x)$. Here, it is sufficient to show that ${d_\omega }(\bar x) = 0$  because of statements (1) and (3) of Lemma 4. Suppose for contradiction that ${d_\omega }(\bar x) \ne 0$. Then, it follows from Lemma 8 that ${\lambda _{{k_j}}} \to 0 $. Therefore, by the definition of ${\lambda _{{k_j}}}$ in Step 4 of Algorithm 1, for sufficiently large $j$ there exists some ${i_{{k_j}}} \in \{ 1,...,m\}$ such that
\[{F_{{i_{{k_j}}}}}({x^{{k_j}}} + {\zeta ^{ - 1}}{\lambda _{{k_j}}}{d^{{k_j}}}) > {F_{{i_{{k_j}}}}}({x^{{k_j}}}) + {\zeta ^{ - 1}}{\lambda _{{k_j}}}\tau {\theta _{{x^{{k_j}}}}}({d^{{k_j}}}).\]

Since $i$ only takes finite number of values in $\{ 1,...,m\} $, we can assume that ${i_{{k_j}}} = \bar i $ without loss of generality. We thus obtain
\begin{align}
\frac{{{F_{\bar i}}({x^{{k_j}}} + {\zeta ^{ - 1}}{\lambda _{{k_j}}}{d^{{k_j}}}) - {F_{\bar i}}({x^{{k_j}}})}}{{{\zeta ^{ - 1}}{\lambda _{{k_j}}}}} > \tau {\theta _{{x^{{k_j}}}}}({d^{{k_j}}}).
\end{align}

Recall that $0 < {\zeta ^{ - 1}}{\lambda _{{k_j}}} < 1$. It follows from the definition   of ${\theta _x}$ that
\begin{align}
{\theta _{{x^{{k_j}}}}}({d^{{k_j}}}) \ge & \nabla {g_{\bar i}}{({x^{{k_j}}})^T}{d^{{k_j}}} + \frac{1}{2}{({d^{{k_j}}})^T}{B_{\bar i}}({x^{{k_j}}})({d^{{k_j}}}) + {h_{\bar i}}({x^{{k_j}}} + {d^{{k_j}}}) - {h_{\bar i}}({x^{{k_j}}})  \notag \\
\ge & \frac{{{\zeta ^{ - 1}}{\lambda _{{k_j}}}\nabla {g_{\bar i}}{{({x^{{k_j}}})}^T}{d^{{k_j}}} + \frac{1}{2}{\zeta ^{ - 1}}{\lambda _{{k_j}}}{{({d^{{k_j}}})}^T}{B_{\bar i}}({x^{{k_j}}})({d^{{k_j}}}) + {h_{\bar i}}({x^{{k_j}}} + {\zeta ^{ - 1}}{\lambda _{{k_j}}}{d^{{k_j}}}) - {h_{\bar i}}({x^{{k_j}}})}}{{{\zeta ^{ - 1}}{\lambda _{{k_j}}}}}  \notag \\
= & \frac{{{g_{\bar i}}({x^{{k_j}}} + {\zeta ^{ - 1}}{\lambda _{{k_j}}}{d^{{k_j}}}) + {h_{\bar i}}({x^{{k_j}}} + {\zeta ^{ - 1}}{\lambda _{{k_j}}}{d^{{k_j}}}) - {g_{\bar i}}({x^{{k_j}}}) - {h_{\bar i}}({x^{{k_j}}}) + o({{({\zeta ^{ - 1}}{\lambda _{{k_j}}}\left\| {{d^{{k_j}}}} \right\|)}^2})}}{{{\zeta ^{ - 1}}{\lambda _{{k_j}}}}}  \notag \\
= & \frac{{{F_{\bar i}}({x^{{k_j}}} + {\zeta ^{ - 1}}{\lambda _{{k_j}}}{d^{{k_j}}}) - {F_{\bar i}}({x^{{k_j}}})}}{{{\zeta ^{ - 1}}{\lambda _{{k_j}}}}} + \frac{{o({{({\zeta ^{ - 1}}{\lambda _{{k_j}}}\left\| {{d^{{k_j}}}} \right\|)}^2})}}{{{\zeta ^{ - 1}}{\lambda _{{k_j}}}}}, \notag
\end{align}
where ${B_{\bar i}}({x^{{k_j}}})$ is some approximation of ${\nabla ^2}{g_{\bar i}}({x^{{k_j}}})$, the second inequality comes from the convexity of ${h_i}$ and Lemma 1, and the first equality follows from the second-order Taylor expansion of ${g_i}$. Therefore, we get
\begin{align}
{\theta _{{x^{{k_j}}}}}({d^{{k_j}}}) \ge \frac{{{F_{\bar i}}({x^{{k_j}}} + {\zeta ^{ - 1}}{\lambda _{{k_j}}}{d^{{k_j}}}) - {F_{\bar i}}({x^{{k_j}}})}}{{{\zeta ^{ - 1}}{\lambda _{{k_j}}}}} + \frac{{o({{({\zeta ^{ - 1}}{\lambda _{{k_j}}}\left\| {{d^{{k_j}}}} \right\|)}^2})}}{{{\zeta ^{ - 1}}{\lambda _{{k_j}}}}}.
\end{align}

From (16) and (17), we get
\[\frac{{{F_{\bar i}}({x^{{k_j}}} + {\zeta ^{ - 1}}{\lambda _{{k_j}}}{d^{{k_j}}}) - {F_{\bar i}}({x^{{k_j}}})}}{{{\zeta ^{ - 1}}{\lambda _{{k_j}}}}} > \tau \frac{{{F_{\bar i}}({x^{{k_j}}} + {\zeta ^{ - 1}}{\lambda _{{k_j}}}{d^{{k_j}}}) - {F_{\bar i}}({x^{{k_j}}})}}{{{\zeta ^{ - 1}}{\lambda _{{k_j}}}}} + \tau \frac{{o({{({\zeta ^{ - 1}}{\lambda _{{k_j}}}\left\| {{d^{{k_j}}}} \right\|)}^2})}}{{{\zeta ^{ - 1}}{\lambda _{{k_j}}}}}.\]

It follows that
\begin{align}
\frac{{{F_{\bar i}}({x^{{k_j}}} + {\zeta ^{ - 1}}{\lambda _{{k_j}}}{d^{{k_j}}}) - {F_{\bar i}}({x^{{k_j}}})}}{{{\zeta ^{ - 1}}{\lambda _{{k_j}}}}} > (\frac{\tau }{{1 - \tau }})\frac{{o({{({\zeta ^{ - 1}}{\lambda _{{k_j}}}\left\| {{d^{{k_j}}}} \right\|)}^2})}}{{{\zeta ^{ - 1}}{\lambda _{{k_j}}}}}.
\end{align}

On the other hand, Lemma 7 yields
\[{\theta _{{x^{{k_j}}}}}({d^{{k_j}}}) \le  - \omega {\left\| {{d^{{k_j}}}} \right\|^2}.\]

Since ${d^{{k_j}}} \to {d_\omega }(\bar x) \ne 0$, it follows from the above inequality and (17) that  there exists   $\gamma = \omega {\left\| d_\omega(\bar x) \right\|^2}   > 0$ such that
\begin{align}
- \gamma  \ge & {\theta _{{x^{{k_j}}}}}({d^{{k_j}}}) \notag \\
\ge & \frac{{{F_{\bar i}}({x^{{k_j}}} + {\zeta ^{ - 1}}{\lambda _{{k_j}}}{d^{{k_j}}}) - {F_{\bar i}}({x^{{k_j}}})}}{{{\zeta ^{ - 1}}{\lambda _{{k_j}}}}} + \frac{{o({{({\zeta ^{ - 1}}{\lambda _{{k_j}}}\left\| {{d^{{k_j}}}} \right\|)}^2})}}{{{\zeta ^{ - 1}}{\lambda _{{k_j}}}}} \notag
\end{align}
for sufficiently large $j$. Therefore, for sufficiently large $j$, the following inequality holds.
\begin{align}
\frac{{{F_{\bar i}}({x^{{k_j}}} + {\zeta ^{ - 1}}{\lambda _{{k_j}}}{d^{{k_j}}}) - {F_{\bar i}}({x^{{k_j}}})}}{{{\zeta ^{ - 1}}{\lambda _{{k_j}}}}} \le  - \gamma  - \frac{{o({{({\zeta ^{ - 1}}{\lambda _{{k_j}}}\left\| {{d^{{k_j}}}} \right\|)}^2})}}{{{\zeta ^{ - 1}}{\lambda _{{k_j}}}}}.
\end{align}

From (18) and (19), we know  that for sufficiently large $j$
\[(\frac{\tau }{{1 - \tau }})\frac{{o({{({\zeta ^{ - 1}}{\lambda _{{k_j}}}\left\| {{d^{{k_j}}}} \right\|)}^2})}}{{{\zeta ^{ - 1}}{\lambda _{{k_j}}}}} <  - \gamma  - \frac{{o({{({\zeta ^{ - 1}}{\lambda _{{k_j}}}\left\| {{d^{{k_j}}}} \right\|)}^2})}}{{{\zeta ^{ - 1}}{\lambda _{{k_j}}}}}.\]

Taking $j \to \infty $, we have $0 <  - \gamma$, which contradicts the fact that $\gamma  > 0$. Therefore, we conclude that ${d_\omega }(\bar x) = 0$.

We now prove the second statement. It is easy to see that the set  $\{ x \in {R^n}\mid F(x) \le F({x^0})\} $ is bounded and that objective values sequence $\{{F}({x^k}) \}$ is $R_+^m$-decrease.  Therefore, the sequence $\{ {x^k}\} $ generated by Algorithm 1 is contained in the above set and so it is also bounded and has at least one accumulation point, which is a Pareto stationary point of (UMOP) according to the first statement.
\end{proof}

\begin{flushleft}
{\zihao{-4}\textbf{4.2~ Convergence of Algorithm 2}}
\end{flushleft}

\textbf{Lemma 9}~~Let $\{ {d^k}\} $ be generated by Algorithm 2 and suppose that $\{ {F_i}({x^k})\} $ is bounded from below for all $i=1,...,m.$ Then, we have
\[\mathop {\lim }\limits_{k \to \infty } {\left\| {{d^k}} \right\|^2} = 0.\]

\begin{proof}[\quad\quad\textbf{Proof}]
From the so-called descent Lemma \cite[Proposition A.24]{B1} and by Lipschitz continuity of $\nabla {g_i}$, we obtain for all $i=1,...,m,$
\[{g_i}({x^k} + {d^k}) \le {g_i}({x^k}) + \nabla {g_i}{({x^k})^T}{d^k} + \frac{L}{2}{\left\| {{d^k}} \right\|^2}.\]
Moreover, since the positive definiteness of ${B_i}({x^k})$ implies ${({d^k})^T}{B_i}({x^k})({d^k}) > 0$ for all $i = 1,...,m.$
Therefore, we get
\begin{align}
{g_i}({x^k} + {d^k}) < {g_i}({x^k}) + \nabla {g_i}{({x^k})^T}{d^k} + \frac{1}{2}{({d^k})^T}{B_i}({x^k})({d^k}) + \frac{L}{2}{\left\| {{d^k}} \right\|^2}.
\end{align}

At the $k$th iteration, we get for $i=1, 2, ..., m$,
\begin{align}
 & {g_i}({x^k} + {d^k}) + {h_i}({x^k} + {d^k}) \notag \\
 = & {g_i}({x^k}) + {h_i}({x^k}) + {g_i}({x^k} + {d^k}) - {g_i}({x^k}) + {h_i}({x^k} + {d^k}) - {h_i}({x^k}) \notag \\
 < & {g_i}({x^k}) + {h_i}({x^k}) + \nabla {g_i}{({x^k})^T}{d^k} + \frac{1}{2}{({d^k})^T}{B_i}({x^k})({d^k}) + {h_i}({x^k} + {d^k}) - {h_i}({x^k}) + \frac{L}{2}{\left\| {{d^k}} \right\|^2} \notag \\
 \le & {g_i}({x^k}) + {h_i}({x^k}) + {\theta _{{x^k}}}({d^k}) + \frac{L}{2}{\left\| {{d^k}} \right\|^2} \notag \\
 \le & {g_i}({x^k}) + {h_i}({x^k}) + \frac{{L - 2\omega }}{2}{\left\| {{d^k}} \right\|^2}. \notag
\end{align}

Here, the first inequality follows from (20), the second inequality  follows from the definition   of ${\theta _x}$, and the third inequality  comes from Lemma 7. Since $\{ {F_i}({x^k})\} $ is bounded from below, there exists ${\hat F_i} \in R$ such that ${\hat F_i} \le {F_i}({x^k}) = {g_i}({x^k}) + {h_i}({x^k})$ for all $i,k.$ Adding up the above inequality from $k = 0$ to $k = \hat k$, where $\hat k$ is a positive integer, we obtain
\[{g_i}({x^{\hat k + 1}}) + {h_i}({x^{\hat k + 1}}) < {g_i}({x^0}) + {h_i}({x^0}) + \frac{{L - 2\omega }}{2}\sum\limits_{k = 0}^{\hat k}{\left\| {{d^k}} \right\|^2}.\]

Since $\omega  > L /2$, we have
\begin{align}
\sum\limits_{k = 0}^{\hat k} {{{\left\| {{d^k}} \right\|}^2}}  < & \frac{2}{{2\omega  - L}}({g_i}({x^0}) + {h_i}({x^0}) - ({g_i}({x^{\hat k + 1}}) + {h_i}({x^{\hat k + 1}})) \notag \\
\le & \frac{2}{{2\omega  - L}}({g_i}({x^0}) + {h_i}({x^0}) - {\hat F_i}). \notag
\end{align}

Taking $\hat k \to \infty $, we obtain
\[\sum\limits_{k = 0}^\infty  {{{\left\| {{d^k}} \right\|}^2}}  < \infty \]

and hence $\mathop {\lim }\limits_{k \to \infty } {\left\| {{d^k}} \right\|^2} = 0 $.
\end{proof}

\textbf{Theorem 2} (i) ~Suppose that $\{ {F_i}({x^k})\} $ is bounded from below for all $i = 1,...,m$.~Then~every accumulation point of the sequence $\{ {x^k}\} $ generated by Algorithm 2 is a Pareto stationary point.

(ii) Moreover, if the level set of  ${F}$ in the sense that $\{ x \in {R^n}\mid F(x) \le F({x^0})\} $ is bounded, then $\{ {x^k}\} $ has accumulation points and they are all Pareto stationary.

\begin{proof}[\quad\quad\textbf{Proof}]
(i) Let $\hat x$ be an accumulation point of $\{ {x^k}\} $ and let $\{{x^{{k_j}}}{\rm{\} }}$ be a subsequence converging to $\hat x$.  From statement (3) of Lemma 4, we have ${d^{{k_j}}} = {d_\omega }({x^{{k_j}}}) \to {d_\omega }(\hat x)$. Here, it is sufficient to show that ${d_\omega }(\hat x) = 0$ because of statements (1) and (3) of Lemma 4. Suppose for contradiction that ${d_\omega }(\hat x) \ne 0$, which contradicts the fact that Lemma 9. Therefore, we conclude that ${d_\omega }(\hat x) = 0$.

(ii) By similar augument with that of (ii) of Theorem 1, it is easy to know the result holds true.
\end{proof}

\vspace{2mm}
\begin{flushleft}
{\zihao{4}\textbf{\Large5\quad Applications}}
\end{flushleft}

\indent
In this section, we consider two applications of the unconstrained multiobjective optimization problem (1) with (2), and discuss how to solve subproblems (8) in a particular application.

\begin{flushleft}
{\zihao{-4}  \textbf{5.1~Application to constrained multiobjective optimization}}
\end{flushleft}

In this section, we consider the following constrained multiobjective optimization problem:
\begin{equation}
\begin{aligned} \label{P}
&\min \quad g(x)  \\
&\begin{array}{r@{\quad}r}
s.t. & x \in S,  \\
\end{array}
\end{aligned}
\end{equation}
where $g:{R^n} \to {R^m}$ is a vector-valued function with $g: = {({g_1},...,{g_m})^T}$ and $S \subset {R^n}$ is convex. Suppose that each component ${g_i}$ of $g$ is twice continuously differentiable. Let $h:{R^n} \to {R^m}$ be a vector-valued function with $h: = {({h_1},...,{h_m})^T}$, where each ${h_i}$ is indicator function of $S$, i.e.,
\[{h_i}(x) = \left\{ \begin{array}{l}
0,x \in S,\\
\infty ,x \notin S.
\end{array} \right.\]

Then, we can rewrite the search direction given in (4) with $x = {x^k}$ as
\[{d^k}: = \mathop {\arg \min }\limits_{d \in S - {x^k}} \left\{ {\max \{ \nabla {g_i}({x^k})d + \frac{1}{2}{d^T}{B_i}({x^k})d\}  + \frac{\omega }{2}{{\left\| d \right\|}^2}} \right\}.\]

\begin{flushleft}
{\zihao{-4}  \textbf{5.2~Application to robust multiobjective optimization}}
\end{flushleft}

Now, let us apply the proposed proximal quasi-Newton algorithms to the robust multiobjective optimization including uncertain parameters
 which is exactly  the  multiobjective optimization problem discussed in \cite{TFY}. In other words, we will solve the (UMOP) with the convex function ${h_i}$ defined as follows:
\begin{align}
{h_i}(x): = \mathop {\max }\limits_{u \in {U_i}} {\hat h_i}(x,u),
\end{align}
where ${U_i} \subseteq {R^n}$ is an uncertainty set, and ${{\hat h}_i}:{R^n} \times {R^n} \to R$ is convex with respect to the first argument.

It is easy to see that ${h_i}$ is also convex. However, ${h_i}$ is not necessarily differentiable even if ${{\hat h}_i}$ is differentiable. First, let us reformulate the subproblem (8) by using an extra variable $\mu  \in R$ as
\begin{equation}
\begin{aligned} \notag
&\mathop {\min }\limits_{\mu ,d} \quad \mu  +  \frac{\omega }{2}{\left\| d \right\|^2}  \\
&\begin{array}{r@{\quad}r}
s.t. &\nabla {g_i}{(x)^T}d + \frac{1}{2}{d^T}{B_i}(x)d + {h_i}(x + d) - {h_i}(x) \le \mu ,i = 1,...,m.  \\
\end{array}
\end{aligned}
\end{equation}

It is worthy noting that ${h_i}$ is not easy to calculate, and thus, the subproblem is difficult to solve.   When ${\hat h_i}$ and ${U_i}$ have some special structure, the constraints of the above problem can be written as explicit formulae by using the duality of (22).  Now, assume that the dual problem of the maximization problem (22) is written as follows:
\begin{equation}
\begin{aligned} \notag
&\mathop {\min }\limits_{{w_i}} \quad {\tilde h_i}(x,{w_i})  \\
&\begin{array}{r@{\quad}r}
s.t. &{w_i} \in {\tilde U_i}(x), \\
\end{array}
\end{aligned}
\end{equation}
where ${{\tilde h}_i}:{R^n} \times {R^m} \to R$ and ${{\tilde U}_i}:{R^n} \to {2^{{R^m}}}$.  If strong duality holds, then we see that the subproblem (8) is equivalent to
\begin{equation}
\begin{aligned} \label{P}
&\mathop {\min }\limits_{\mu ,d,{w_i}} \quad \mu  +  \frac{\omega }{2}{\left\| d \right\|^2} \\
&\begin{array}{r@{\quad}r}
s.t. &\nabla {g_i}{(x)^T}d + \frac{1}{2}{d^T}{B_i}(x)d + {\tilde h_i}(x + d,{w_i}) - {h_i}(x) \le \mu , \\
     &{w_i} \in {\tilde U_i}(x + d),i = 1,...,m. \\
\end{array}
\end{aligned}
\end{equation}

When ${\tilde h_i}$ and ${\tilde U_i}$ have some explicit form, this problem is tractable. As we mention below, in this case, we can convert the above subproblem to some well-known convex optimization problems. This idea can be also seen in \cite{TFY}. In the following, we will introduce a robust multiobjective optimization problem where the subproblem can be written as a quadratic programming.

Suppose that ${{\hat h}_i}(x,u) = {u^T}x$ and ${U_i} = \{ u \in {R^n}|{A_i}u \le {b_i}\} $, where ${A_i} \in {R^{d \times n}}$ and ${b_i} \in {R^d}$, that is, ${{\hat h}_i}$ is linear in $x$, and ${U_i}$ is a polyhedron. Suppose also that ${U_i}$ is nonempty and bounded. Then, follow the ideas of  Tanabe, Fukuda and  Yamashita \cite{TFY}, problem (22) can be  rewritten as  the following linear programming problem:
\begin{equation}
\begin{aligned} \label{P}
&\mathop {\max }\limits_u \quad {x^T}u  \\
&\begin{array}{r@{\quad}r}
s.t. &{A_i}u \le {b_i}, \\
\end{array}
\end{aligned}
\end{equation}
and its dual problem is given by
\begin{equation}
\begin{aligned} \notag
&\mathop {\min }\limits_w \quad {b_i}^Tw  \\
&\begin{array}{r@{\quad}r}
s.t. &{A_i}^Tw = x , \\
     &w \ge 0. \\
\end{array}
\end{aligned}
\end{equation}

Since the strong duality holds, we can convert the subproblem (8) [or, equivalently (23)] to a nonlinearly constrained quadratic programming problem:
\begin{equation}
\begin{aligned} \label{P}
&\mathop {\min }\limits_{\mu ,d,{w_i}} \quad \mu  + \frac{\omega }{2}{\left\| d \right\|^2} \\
&\begin{array}{r@{\quad}r@{}}
s.t. &\nabla {g_i}{(x)^T}d + \frac{1}{2}{d^T}{B_i}(x)d + {b_i}^T{w_i} - {h_i}(x) \le \mu , \\
     &{A_i}^T{w_i} = x + d, \\
     &{w_i} \ge 0,i = 1,...,m. \\
\end{array}
\end{aligned}
\end{equation}

\vspace{2mm}
\begin{flushleft}
{\zihao{4}\textbf{\Large6\quad Numerical experiments}}
\end{flushleft}

In this section, we present some numerical results using Algorithms 1 and 2 for the nonlinearly constrained quadratic programming problem in Section 5.2. The experiments are carried out on a machine with a 2.2GHz Intel Core i3 CPU and 6GB memory, and we implement all codes in MATLAB R2018b. We
consider the problem (UMOP) (i.e., the problem (1) with (2)), where $n = 5$, $m = 2$, ${g_i}(x) = \frac{1}{2}{x^T}{Q_i}x + {q_i}^Tx$, ${h_i}(x) = \mathop {\max }\limits_{u \in {U_i}} {\hat h_i}(x,u)$, ${Q_i} \in {R^{n \times n}}$, ${q_i} \in {R^n}$, and  ${{\hat h}_i}:{R^n} \times {R^n} \to R ,i = 1,...,m$. Here, we assume that each ${Q_i}$ is positive definite, so it can be decomposed as ${Q_i} = {M_i}{M_i}^T$, where ${M_i} \in {R^{n \times n}}$. We generate ${M_i}$ and ${q_i}$ by choosing every component randomly from the standard normal distribution. To implement Algorithms 1 and 2, we make the following choices.

\textbf{Remark 2}~~(1) Every component of ${x^0}$ is chosen randomly from the standard normal distribution.

                   (2) $B_i^0$ is the identity matrix, i.e., $B_i^0 = I,i = 1,...,m$.

                   (3) We set the constant $\omega  = 5$, $\tau  = \frac{1}{2}$, and $\varsigma  = \frac{1}{2}$.

                   (4) The terminate criteria is replaced by $\left\| {{d^k}} \right\| < \varepsilon : = {10^{ - 6}}$.

Also, we run each one of the following experiments 100 times from different initial points, and with $\delta  = 0,~ 0.05,~ 0.1$. Naturally, when $\delta  = 0$, no uncertainties are considered.

In order to solve the nonlinear constrained quadratic programming problem (25), the following numerical experiments are performed. We assume that ${h_i}(x) = \mathop {\max }\limits_{u \in {U_i}} {u^T}x,i = 1,2$, where ${U_1} = \{ u \in {R^5}| - \delta  \le {u_i} \le \delta ,i = 1,...,5\} $ and ${U_2} = \{ u \in {R^5}| - \delta  \le {(Bu)_i} \le \delta ,i = 1,...,5\} $. Here, every component of $B \in {R^{5 \times 5}}$ is chosen randomly from the standard normal distribution and $\delta  \ge 0$.

For simplicity of the notation, from now on, the proximal gradient method for (UMOP), which was introduced in \cite{TFY} is denoted by PGM,  the proximal quasi-Newton methods for (UMOP) related to   BFGS, self-scaling BFGS, and Huang BFGS  methods   are denoted by PQNM(BFGS), PQNM(SS-BFGS), and PQNM(H-BFGS), respectively. We will use the toolbox of convex optimization in MATLAB to solve (24) and (25).

\begin{flushleft}
{\zihao{-4}\textbf{6.1~ Numerical experiment of algorithm 1}}
\end{flushleft}

On one hand, the experimental result obtained by using PGM with line searches (i.e., Algorithm 3.1 in \cite{TFY}) with  different $\delta $ for (UMOP) is shown in Fig. 1. The experimental results obtained by using PQNM(BFGS), PQNM(SS-BFGS), and PQNM(H-BFGS) with line searches for (UMOP) with deferent $\delta $ is shown in Fig. 2, Fig. 3 and Fig. 4, respectively. For each $\delta $, we obtained part of the Pareto frontier of (UMOP). And it is easy to see that    the Pareto frontier of (UMOP) becomes lower when $\delta $ is smaller. On the other hand, Fig. 5 shows the comparison among the (PGM), PQNM(BFGS), PQNM(SS-BFGS) and PQNM(H-BFGS) with line searches for (UMOP) when $\delta $ is 0, 0.05 and 0.1, respectively. When $\delta $ is same, we observed that  the Pareto frontier of (UMOP) produced by PQNM(H-BFGS) with line searches  is the lowest, the Pareto frontier of (UMOP) produced by  PGM with line searches   is the highest, and the Pareto frontier of (UMOP) produced by   PQNM(SS-BFGS) with line searches   is lower than that produced by the PQNM(BFGS) with line searches. And so  PQNM(H-BFGS) with line searches for (UMOP) is the most effective method and the introduced new proximal quasi-Newton methods with line searches for (UMOP) are more effective than PGM with line searches for (UMOP).

\begin{figure}[H]

  \subfigure{
  \begin{minipage}{7cm}
    \centering
      \includegraphics[height=4.8cm,width=7.8cm]{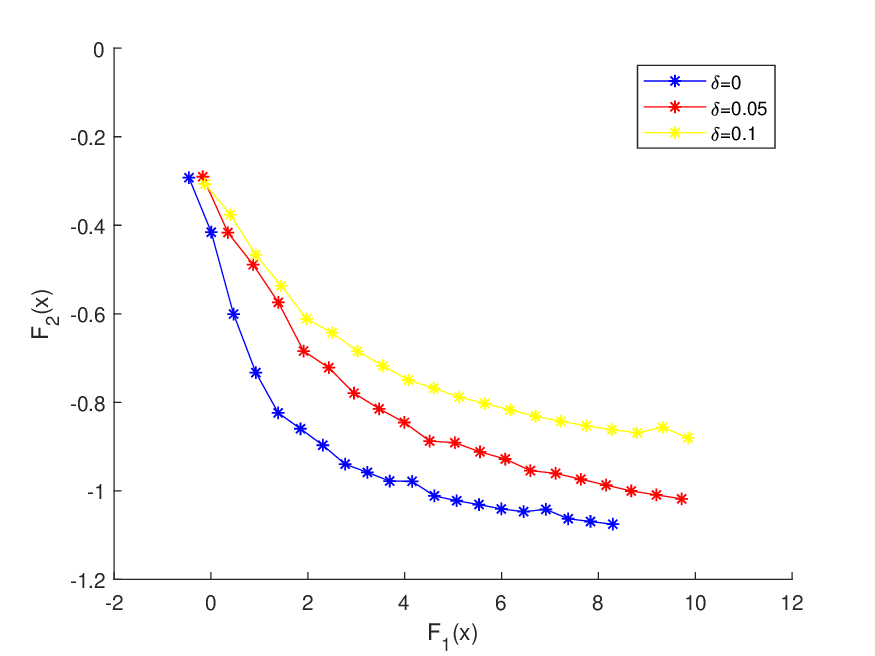}
      \begin{center}
      \textbf{Fig.1}~~Result for PGM
      \end{center}
    \end{minipage}%
  }
  \subfigure{
    \begin{minipage}{7cm}
        \centering
      \includegraphics[height=4.8cm,width=7.8cm]{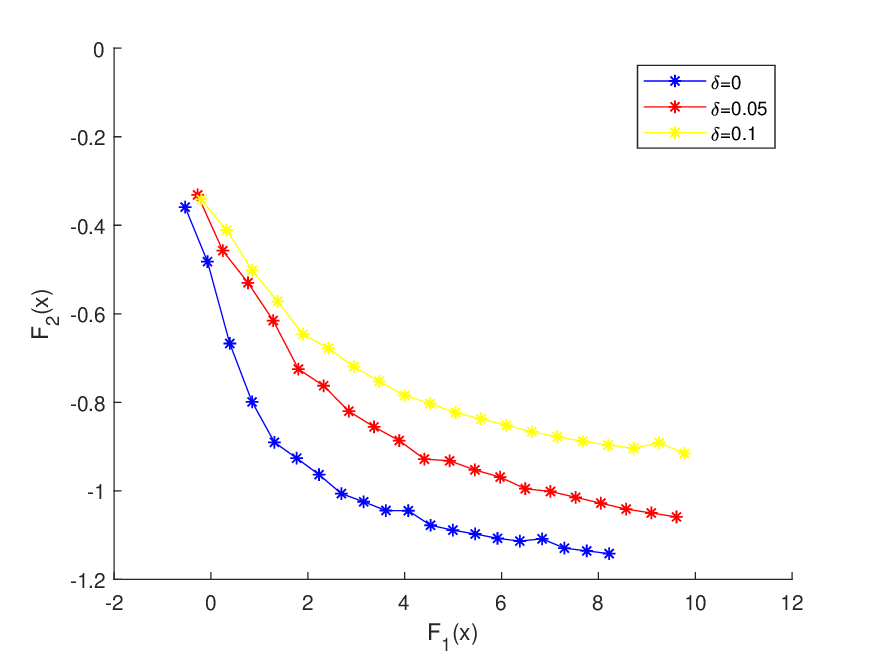}
      \begin{center}
      \textbf{Fig.2}~~Result for PQNM(BFGS)
      \end{center}
    \end{minipage}%
  }
\end{figure}

\begin{figure}[H]

  \subfigure{
  \begin{minipage}{7cm}
    \centering
      \includegraphics[height=4.8cm,width=7.8cm]{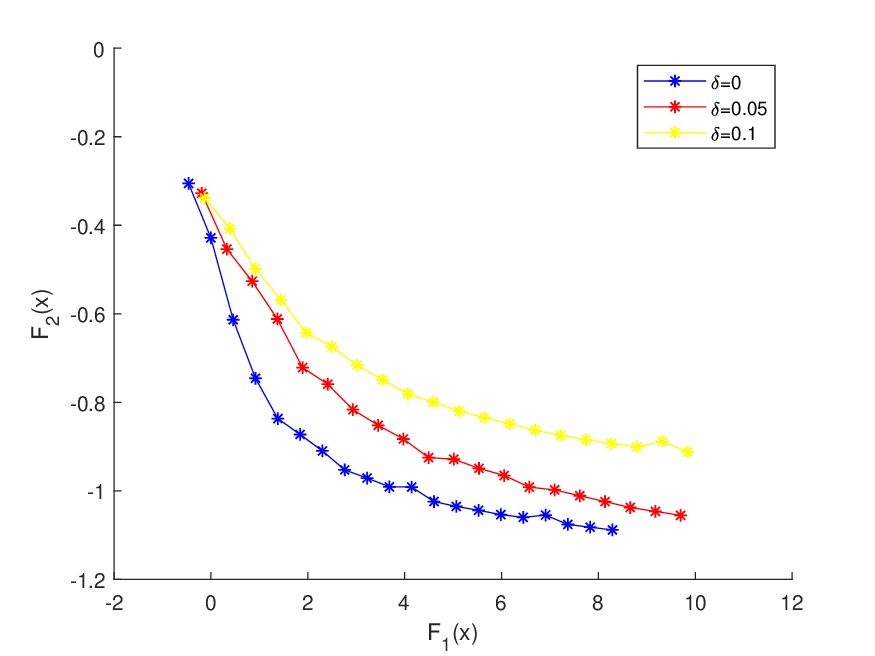}
      \begin{center}
      \textbf{Fig.3}~~Result for PQNM(SS-BFGS)
      \end{center}
    \end{minipage}%
  }
  \subfigure{
    \begin{minipage}{7cm}
        \centering
      \includegraphics[height=4.8cm,width=7.8cm]{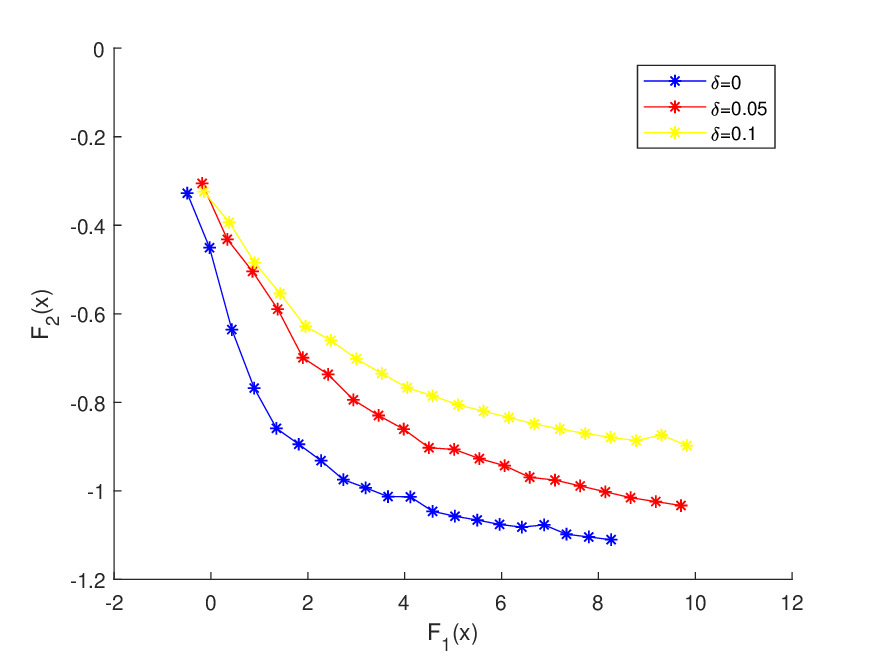}
      \begin{center}
      \textbf{Fig.4}~~Result for PQNM(H-BFGS)
      \end{center}
    \end{minipage}%
  }
\end{figure}

\begin{flushleft}
{\zihao{-4}\textbf{6.2~Numerical experiment of algorithm 2}}
\end{flushleft}

Firstly, the experimental result of the PGM without line searches (i.e., Algorithm 3.2 in \cite{TFY}) for (UMOP) is shown in Fig. 6. And  the experimental results obtained by using PQNM(BFGS), PQNM(SS-BFGS), and PQNM(H-BFGS) without line searches for (UMOP) with deferent $\delta $ is shown in Fig. 7, Fig. 8, Fig. 9, respectively. It is obvious that the Pareto frontier produced by the same method become smaller as $\delta $ gets smaller. Secondly, Fig.10 shows the comparison among the (PGM), PQNM(BFGS), PQNM(SS-BFGS) and PQNM(H-BFGS) without line searches for (UMOP) when $\delta $ is 0, 0.05 and 0.1, respectively. When $\delta $ is  same, we also observed that the Pareto frontier of (UMOP) produced by PQNM(H-BFGS) without line searches  is the lowest one, the Pareto frontier of   (UMOP) produced by PGM without line searches for (UMOP) is the highest one, and the Pareto frontier of   (UMOP) produced by  PQNM(SS-BFGS) without line searches   is lower than that produced by the    PQNM(BFGS) without line searches. And so PQNM(H-BFGS) without line searches for (UMOP) is the most effective method and the introduced new proximal quasi-Newton methods without line searches for (UMOP) are more effective than PGM without line searches for (UMOP).

\begin{figure}[H]
  \centering
  \subfigure{
  \begin{minipage}{7cm}
    \centering
      \includegraphics[height=4.8cm,width=7.8cm]{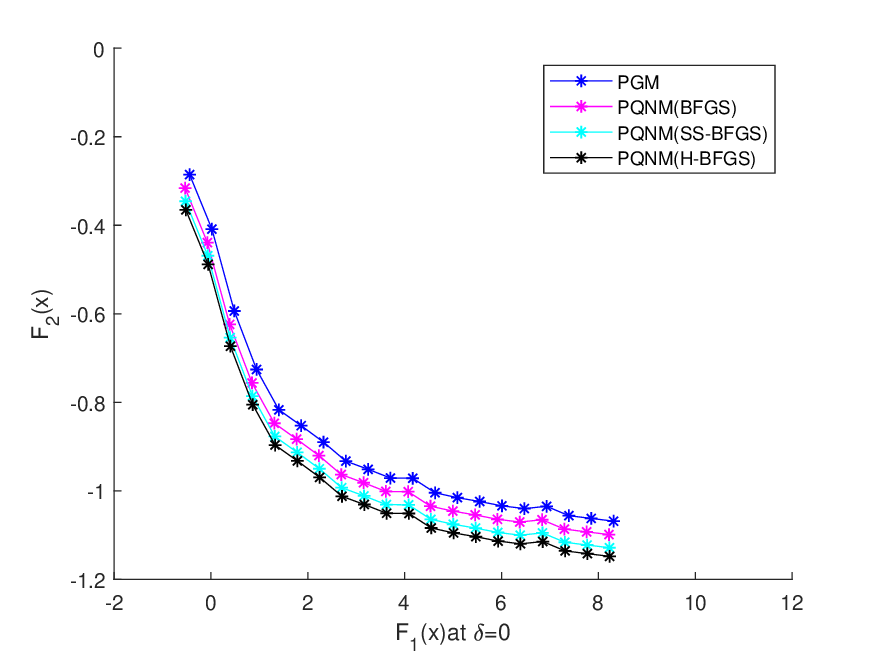}
    \end{minipage}%
  }
  \subfigure{
    \begin{minipage}{7cm}
        \centering
      \includegraphics[height=4.8cm,width=7.8cm]{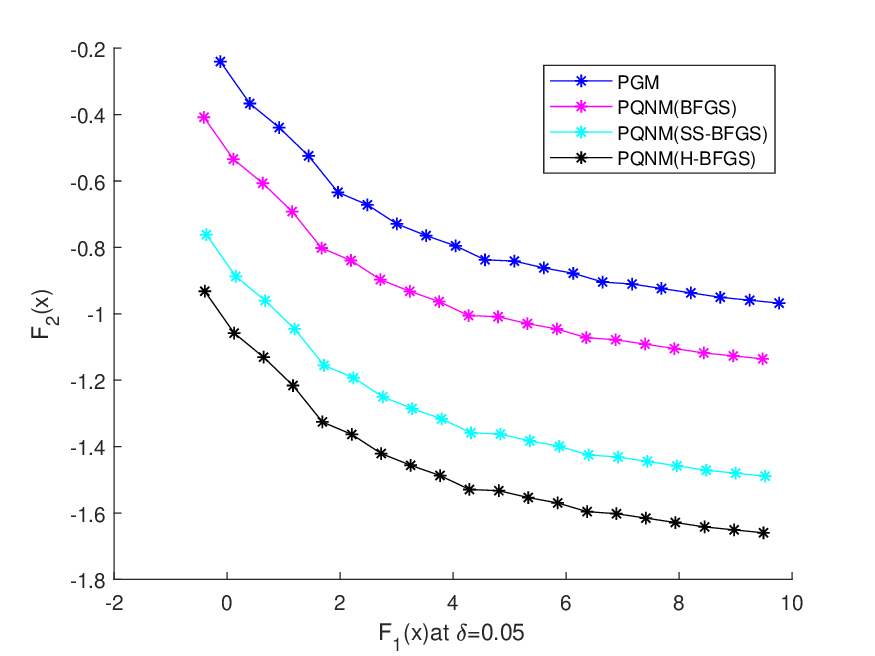}
    \end{minipage}%
  }
  \subfigure{
    \begin{minipage}{7cm}
     \centering
      \includegraphics[height=4.8cm,width=7.8cm]{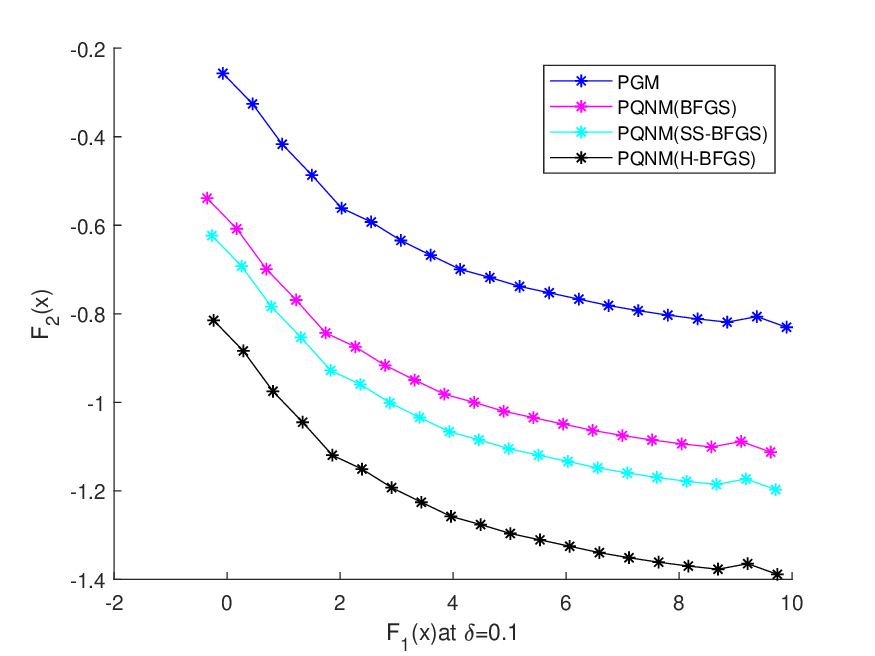}
    \end{minipage}%
  }
   \begin{center}
      \textbf{Fig.5}~~Comparison of PGM, PQNM(BFGS), PQNM(SS-BFGS) and PQNM(H-BFGS) when $\delta$ is 0, 0.05 and 0.1, respectively
   \end{center}
\end{figure}

\begin{figure}[H]

  \subfigure{
  \begin{minipage}{7cm}
    \centering
      \includegraphics[height=4.8cm,width=7.8cm]{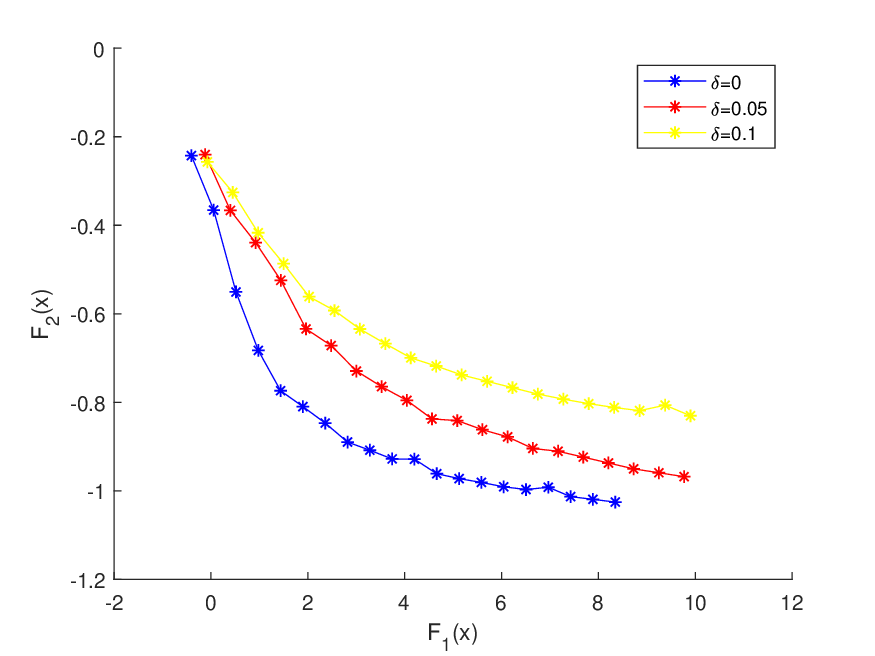}
      \begin{center}
      \textbf{Fig.6}~~Result for PGM
      \end{center}
    \end{minipage}%
  }
  \subfigure{
    \begin{minipage}{7cm}
        \centering
      \includegraphics[height=4.8cm,width=7.8cm]{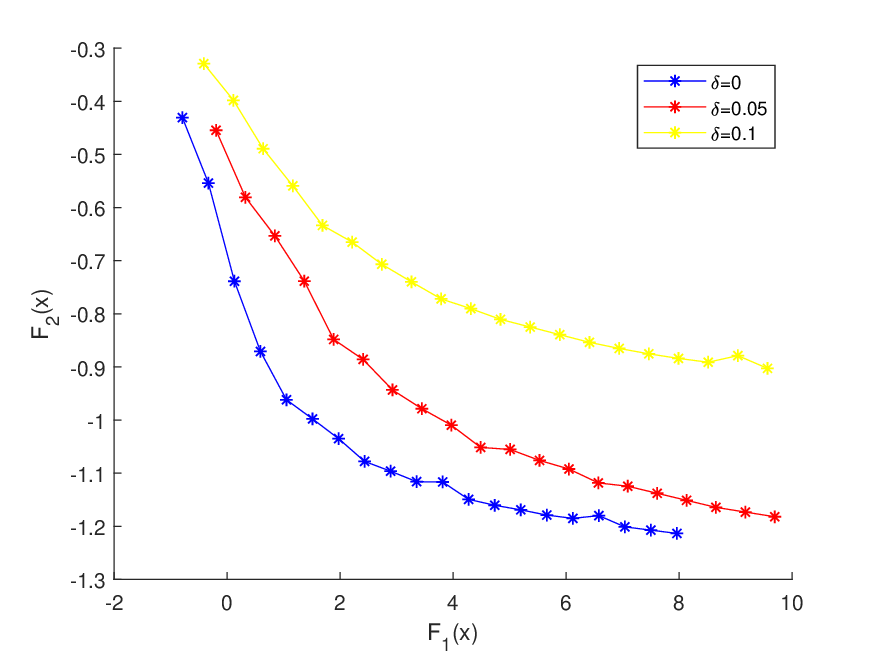}
      \begin{center}
      \textbf{Fig.7}~~Result for PQNM(BFGS)
      \end{center}
    \end{minipage}%
  }
\end{figure}

\begin{figure}[H]

  \subfigure{
  \begin{minipage}{7cm}
    \centering
      \includegraphics[height=4.8cm,width=7.8cm]{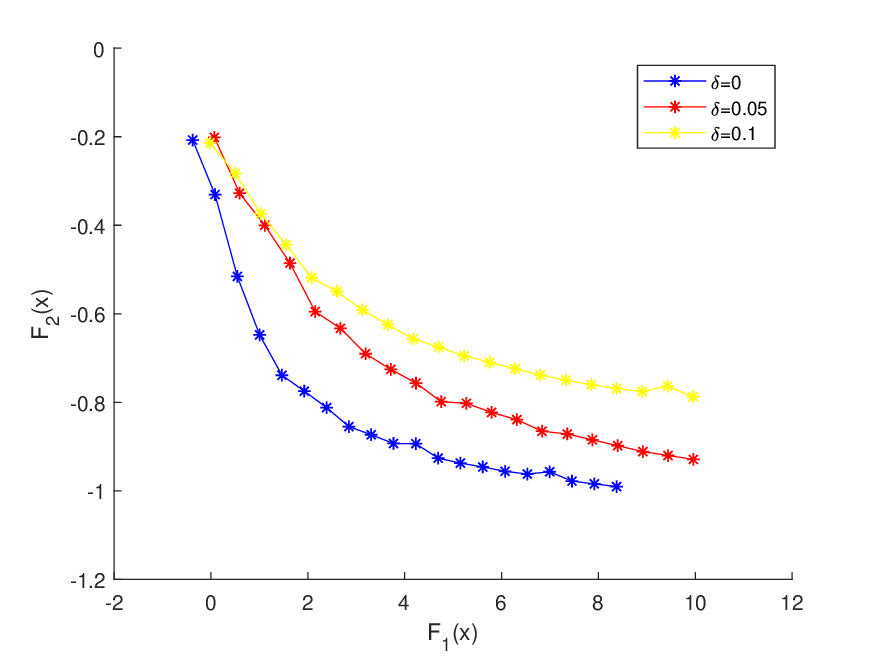}
      \begin{center}
      \textbf{Fig.8}~~Result for PQNM(SS-BFGS)
      \end{center}
    \end{minipage}%
  }
  \subfigure{
    \begin{minipage}{7cm}
        \centering
      \includegraphics[height=4.8cm,width=7.8cm]{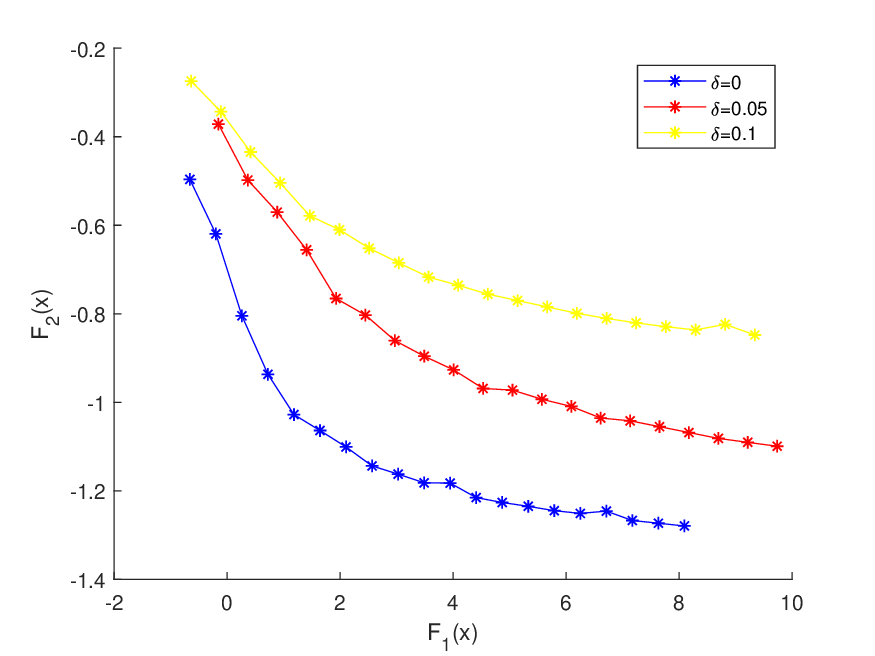}
      \begin{center}
      \textbf{Fig.9}~~Result for PQNM(H-BFGS)
      \end{center}
    \end{minipage}%
  }
\end{figure}

\begin{figure}[H]
  \centering
  \subfigure{
  \begin{minipage}{7cm}
    \centering
      \includegraphics[height=4.8cm,width=7.8cm]{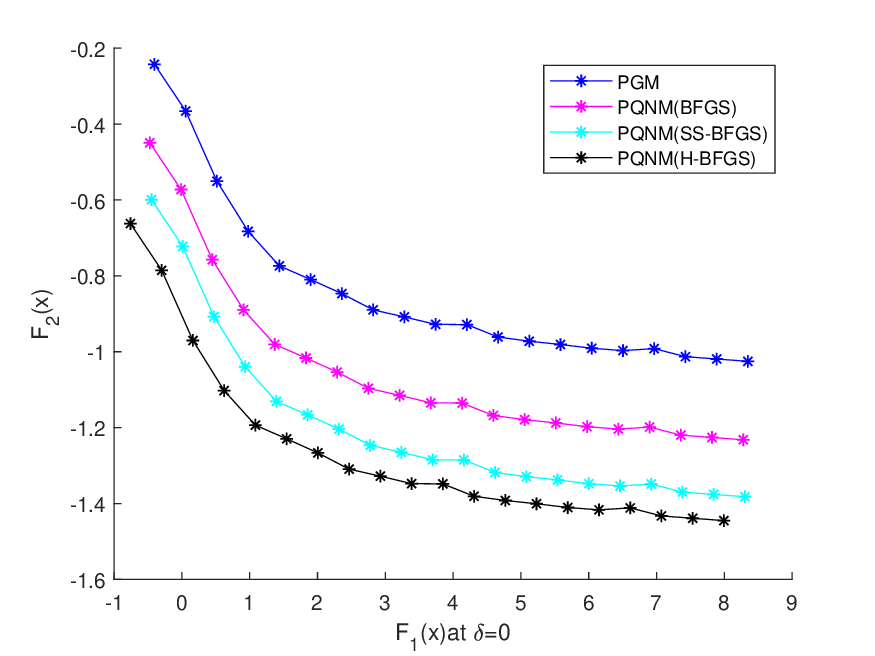}
    \end{minipage}%
  }
  \subfigure{
    \begin{minipage}{7cm}
        \centering
      \includegraphics[height=4.8cm,width=7.8cm]{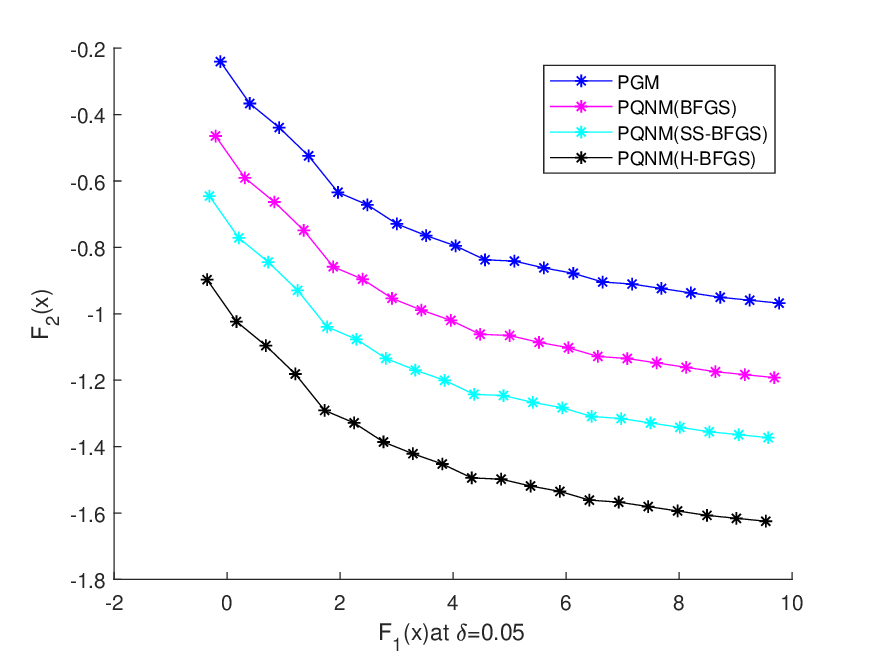}
    \end{minipage}%
  }
  \subfigure{
    \begin{minipage}{7cm}
     \centering
      \includegraphics[height=4.8cm,width=7.8cm]{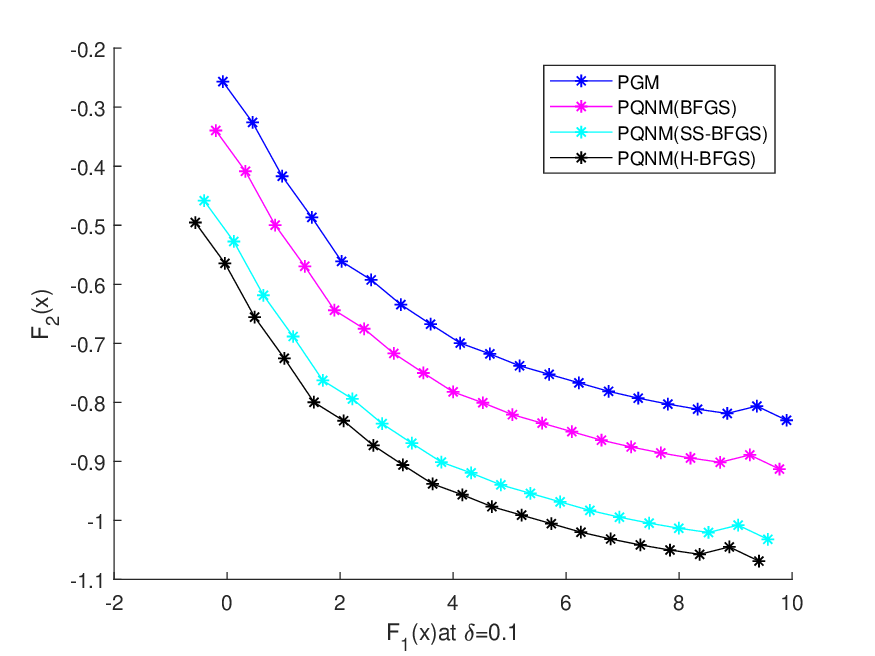}
    \end{minipage}%
  }
   \begin{center}
      \textbf{Fig.10}~~Comparison of PGM, PQNM(BFGS), PQNM(SS-BFGS) and PQNM(H-BFGS) when $\delta$ is 0, 0.05 and 0.1, respectively
   \end{center}
\end{figure}

\vspace{2mm}
\begin{flushleft}
{\zihao{4}\textbf{\Large7\quad Conclusion}}
\end{flushleft}

First, for unconstrained multiobjective optimization problems, where each objective function is the sum of a twice continuously differentiable strongly convex function and a proper convex but not necessarily differentiable function, the proximal quasi-Newton method with line searches and the proximal quasi-Newton method without line searches are proposed. Secondly, under appropriate conditions, we prove that each cluster point of the sequence generated by these two types of algorithms is the Pareto stationary point of the multiobjective optimization problem. Thirdly, we give their applications in constrained multiobjective optimization and robust multiobjective optimization. In particular, for robust multiobjective optimization, we show that the subproblems of the proximal quasi-Newton method can be regarded as quadratic programming problems. Finally, numerical experiments are performed to verify the effectiveness of the proposed algorithms. In the future, we will analyze the convergence rate of the proposed algorithms.

\renewcommand\refname


{\flushleft References}
\begin{thebibliography}{99}
\makeatletter
\renewcommand\@biblabel[1]{#1.}
\makeatother
\bibitem{BLM} Bello Cruz, J.Y., Lucambio Pérez, L.R., Melo, J.G.: Convergence of the projected gradient method for quasiconvex multiobjective optimization. Nonlinear Anal. {\bf 74}(16), 5268--5273 (2011)
\bibitem{BP} Berge, C., Patterson, E.M.: Topological Spaces. Dover Publications, Edinburgh (1963)
\bibitem{B1} Bertsekas, D.P.: Nonlinear Programming, 2nd edn. Athena Scientific, Belmont (1999)
\bibitem{BIS} Bonnel, H., Iusem, A.N., Svaiter, B.F.: Proximal methods in vector optimization. SIAM J. Optim. {\bf 15}(4), 953--970 (2005)
\bibitem{B2} Broyden, C.G.: A new double-rank minimization algorithm. Not. Am. Math. Soc. {\bf 16}, 670 (1969)
\bibitem{CT} Chen, G., Teboulle, M.: Convergence analysis of a proximal-like minimization algorithm using Bregman functions. SIAM J. Optim. {\bf 3}(3), 538--543 (1993)
\bibitem{DDFL} Da Cruz Neto, J.X., Da Silva, G.J.P., Ferreira, O.P., Lopes, J.O.: A subgradient method for multiobjective optimization. Comput. Optim. Appl. {\bf 54}(3), 461--472 (2013)
\bibitem{F} Fletcher, R.: A new approach to variable metric algorithms. Comput. J. {\bf 13}(3), 317--322 (1970)
\bibitem{FGS} Fliege, J., Grana Drummond, L.M., Svaiter, B.F.: Newton’s method for multiobjective optimization. SIAM J. Optim. {\bf 20}(2), 602--626 (2009)
\bibitem{FS} Fliege, J., Svaiter, B.F.: Steepest descent methods for multicriteria optimization. Math. Methods Oper. Res. {\bf 51}(3), 479--494 (2000)
\bibitem{FG} Fukuda, E.H., Grana Drummond, L.M.: A survey on multiobjective descent methods. Pesquisa Operacional {\bf 34}(3), 585--620 (2014)
\bibitem{GS} Gass, S., Saaty, T.: The computational algorithm for the parametric objective function. Naval Res. Logist. Q. {\bf 2}(1--2), 39--45 (1955)
\bibitem{G1} Geoffrion, A.M.: Proper efficiency and the theory of vector maximization. J. Math. Anal. Appl. {\bf 22}(3), 618--630 (1968)
\bibitem{G2} Goldfarb, D.: A family of variable-metric methods derived by variational means. Math. Comp. {\bf 24}, 23--26 (1970)
\bibitem{H} Hogan, W.W.: Point-to-set maps in mathematical programming. SIAM Rev. {\bf 15}(3), 591--603 (1973)
\bibitem{NW} Nocedal, J., Wright, S.J.: Numerical Optimization, 2nd ed. Springer Science and Business Media, LLC, New York (2006)
\bibitem{P} Povalej, Z.: Quasi-Newton's method for multiobjective optimization. J. Comput. Appl. Math. {\bf 255}, 765--777 (2014)
\bibitem{R} Rockafellar, R.T.: Convex Analysis. University Press, Princeton (1970)
\bibitem{S} Shanno, D.F.: Conditioning of quasi-Newton methods for function minimization. Math. Comp. {\bf 24}, 647--656 (1970)
\bibitem{SY} Sun, W., Yuan, Y.X.: Optimization theory and methods: nonlinear programming. Springer, New York (2006)
\bibitem{TFY} Tanabe, H., Fukuda, E.H., Yamashita, N.: Proximal gradient methods for multiobjective optimization and their applications. Comput. Optim. Appl. {\bf 72}, 339--361 (2019)
\bibitem{Z} Zadeh, L.: Optimality and non-scalar-valued performance criteria. IEEE Trans. Autom. Control {\bf 8}(1), 59--60 (1963)
\bibitem{ZX} Zhang, J., Xu, C.: Properties and numerical performance of quasi-Newton methods with modified quasi-Newton equations. J. Comput. Appl. Math. {\bf 137}(2), 269--278 (2001)
\end{thebibliography}
\end{document}